\documentclass[dvips, 10pt, a4paper]{article}
\usepackage{amsmath,amsfonts,amssymb,amsthm,graphics,amscd}

\font\bigbf=cmbx12

\begin{document}

\centerline{\bigbf SLANT SUBMANIFOLDS OF LORENTZIAN}
 \smallskip
\centerline{\bigbf ALMOST CONTACT MANIFOLDS}
\smallskip

 \bigskip

\centerline{\bf  Khushwant Singh, Siraj Uddin, Satvinder Singh
Bhatia and M. A. Khan}

\thispagestyle{empty} 

\footnote[0]{2000 {\it Mathematics Subject Classifications}.
53C40, 53B25, 53C15.} \footnote[0]{{\it Key words and Phrases}. Space like submanifold, slant submanifold, Lorentzian almost contact manifold, Lorentzian Sasakian.}


\begin{abstract}
In this paper we study slant submanifolds of Lorentzian almost contact manifolds. We have taken the submanifold as a space like and then defined the slant angle on a submanifold and thus we extended the results of [7] and [8] in this new setting.\\
\end{abstract}

\section{Introduction}

\parindent=8mm
Slant submanifolds were introduced by B.Y. Chen in [5, 6]. These are
generalization of both holomorphic and totally real submanifolds of
an almost Hermitian manifold. Since then many research articles have
appeared on these submanifolds in different known spaces. A. Lotta
[8] defined and studied slant submanifolds in contact geometry.
Later on, J.L. Cabrerizo, A. Carriazo, L.M. Fernandez and M.
Fernandez studied slant submanifolds of Sasakian manifolds [4].
Recently, M.A. Khan et.al [7] studied these submanifolds in the
setting of Lorentzian paracontact manifolds.

\parindent=8mm
In this paper we defined and studied slant submanifolds of Lorentzian almost contact manifolds. In section 2, we review some formulae for Lorentzian almost contact manifolds and their submanifolds. In section 3, we define slant submanifold assuming that it is space like for Lorentzian almost contact manifolds. we have given in this section characterization theorems for slant submanifold in the setting of Lorentzian almost contact manifolds. The section 4, has been devoted to the study of slant submanifolds of Lorentzian Sasakian manifolds.

\section{Preliminaries}

\parindent=8mm
Let $\bar M$ be a $(2n+1)-$dimensional manifold with an almost
contact structure and compatible Lorentzian metric, $(\bar M,\phi,
\xi, \eta, g)$ that is, $\phi$ is $(1,1)$ tensor field , $\xi$ is a
structure vector field, $\eta$ is $1$-form and $g$ is Lorentzian
metric on $\bar M$ satisfying [1,3]
$$\phi^2 X=-X+\eta (X)\xi,~~\eta(\xi)=1,~~\phi(\xi)=0,~~\eta\circ\phi=0\eqno(2.1)$$
and
$$g(\phi X, \phi Y)=g(X, Y)+\eta(X)\eta(Y),~~\eta(X)=-g(X, \xi)\eqno(2.2)$$
for any $X, Y\in T\bar M$, where $T\bar M$ denotes the Lie algebra of vector fields on $\bar M$. An almost contact manifold with Lorentzian metric $g$ is called a {\it{Lorentzian almost contact manifold}}. From (2.2), it follows that
$$ g(\phi X, Y)=-g(X, \phi Y).\eqno(2.3)$$
A Lorentzian almost contact manifold is {\it{Lorentzian Sasakian}} if [1]
$$(\bar\nabla_X\phi)Y=-g(X, Y)\xi-\eta(Y)X.\eqno(2.4)$$
It is easy to compute from (2.4) that
$$\bar\nabla_X\xi=\phi X.\eqno(2.5)$$

\parindent=8mm
Now, let $M$ be a submanifold of $\bar M$, we denote the induced Lorentzian metric on $M$ by the same symbol $g$. Let $\bar\nabla$ and $\nabla$ be Levi-Civita connections on the ambient manifold $\bar M$ and the submanifold $M$, respectively with respect to the Lorentzian metric $g$ then the Gauss and Weingarten formulae are given by
$$\bar\nabla_XY=\nabla_X Y+h(X,Y) \eqno(2.6)$$
$$\bar \nabla_XV=-A_VX+\nabla_X^\perp V\eqno(2.7)$$
for any $X,Y \in TM$ and $V\in T^\perp M$, where $\nabla^\perp$ is
the connection on the normal bundle $T^\perp M$, $h$ is the second
fundamental form and $A_V$ is the Weingarten map associated with $V$
as [10]
$$g(A_V X,Y)=g(h(X,Y),V). \eqno(2.8)$$
for any $x\in M$, $X\in T_x M$ and $V\in T^\perp_x M$, we write
$$\phi X=TX+NX\eqno(2.9)$$
$$\phi V=tV+nV\eqno(2.10)$$
where $TX$(resp. $tV$) denotes the tangential component of
$\phi X$(resp. $\phi V$) and $NX$(resp. $nV$) denotes the normal
component of $\phi X$(resp.$\phi V$).

\section{Slant submanifolds}

\noindent
Throughout, this section we consider a submanifold
$M$ of a Lorentzian manifold $\bar M$ such that for all $X\in TM$,
$g(X,X)>0$ or $g(X,X)=0$ i.e., all the tangent vectors on $M$ are
Space like or null like, we shall call these type of submanifolds as {\it{space like}} and also we assume that the structure vector
field $\xi$ is tangent to the submanifold $M$. Let $M$ be an immersed
submanifold of $\bar M$ and for any $x\in M$ and $X\in T_xM$, if the
vector field $X$ and $\xi$ are linearly independent then the angle
$\theta(X)\in [0,\pi/2]$ between $\phi X$ and $T_xM$ is well
defined, if $\theta(X)$ does not depend upon the choice of $x\in M$
and $X\in T_xM$, then $M$ is {\it{slant}} in $\bar M$. The constant angle
$\theta(X)$ is then called the {\it{slant angle}} of $M$ in $\bar M$ and
which in short we denote by Sla$(M)$. The
tangent bundle $TM$ at every point $x\in M$ is decomposed as
$$TM=D\oplus\langle\xi\rangle$$
where $\langle\xi\rangle$ is the one dimensional distribution orthogonal to the slant distribution $D$ on $M$ and spanned by the structure vector field $\xi$.

\parindent=8mm
For any $x\in M$ taking $X\in T_X M$ we put $\phi X=TX+NX$ where $TX\in T_X M$ and $NX\in T^\perp_x M$. Thus defining an endomorphism $T:T_xM\longrightarrow T_xM $, whose square $T^2$
will be denoted by $Q$. Then tensor fields on $M$ of the type $(1,1)$
determined by their endomorphisms shall be denoted by same letters $T$
and $Q$. It is easy to show that for every $x\in M$,
$g(TX,Y)=-g(X,TY)$, which implies that $Q$ is anti-symmetric. Moreover,
in the following steps we can prove that the eigen value of $Q$
always belong to $[-1,0]$. For any $X\in T_x M-<\xi>$, we get
$$g(QX,X)=-\|TX\|^2$$
but,
$$\|TX\|\leq\|QX\|$$
$$\|TX\|\leq \mu \|X\|,~~{\mbox{and}}~~\mu\in [0,1].$$
Thus we obtain
$$g(QX,X)=-\mu(X)\|X\|^2.$$
That is,
$$g(QX,X)=\lambda(X)\|X\|^2$$
where $-1\leq\lambda(X)\leq0$ and $\lambda$ depends on
$X$. In other words, each eigen value of $Q$ lies in $[-1,0]$ and
each eigen value has even multiplicity.

\parindent=8mm Now, we have the following theorem.\\

\noindent{\bf{Theorem 3.1.}} {\it{Let $x\in M$ and $X\in T_x M$ be an
eigenvector of $Q$ with eigenvalue $\lambda(X)$. Suppose $X$ is
linearly independent from $\xi_x$, then,}}
$$\cos\theta(X)=\sqrt{-\lambda(X)}\frac{\|X\|}{\|\phi X\|}.\eqno(3.1)$$\\

\noindent{\it{Proof.}} For any $X\in TM$ we have
$$\|TX\|^2= g(TX,TX)=-\lambda(X)\|X \|^2.\eqno(3.2)$$
On the other hand by definition of $\theta (X)$ we have
$$\cos\theta(X)=\frac {g(\phi X,TX)}{\|TX\|\|\phi X\|}~~~~~~~~~~~~~~~~~~~~~~~~~~~$$
$$~~~~~~~~~~~~~~=\frac{g( TX, TX)}{\|TX\|\|\phi X\|}=-\lambda(X)\frac{\|X\|^2}{\|TX\|\|\phi X\|}.$$
Then from (3.2), we obtain that
$$\cos(\theta)(X)=\sqrt{-\lambda(X)}\frac{\|X\|}{\|\phi X\|}.$$
This completes the proof.$~~~~~~~~~~~~~~~~~~~~~~~~~~~~~~~~~~~~~~~~~~~~~~~~~~~~~~~~~~~~~~~~~~~~~~~~~~~~~~~~~~~~~~~~~~\square$

\parindent=8mm The following characterization theorem gives the existence of eigen values of the endomorphism $Q$.\\

\noindent{\bf{Theorem 3.2.}} {\it{Let $M$ be a space
like slant submanifold of a Lorentzian almost contact manifold $\bar
M$ and $\theta= Sla(M) \neq \pi/2$, then $Q$ admits the real number
$-\cos^2 \theta$ as the only non-vanishing eigen value, for any
$x\in M$. Moreover the related eigen space $H$ satisfies $H\subset
D$, where $D=Span(\xi_x)^\perp \subset T_x M $.}}\\

\noindent {\it{Proof.}} Let $x \in M$, from equation (3.1) $Ker
(Q)\neq T_x M$, otherwise $Sla(M)=\pi/2$ which contradict the
assumption. So let $\lambda$ be an arbitrary non-vanishing eigen
value of $Q$ and let $H$ be the corresponding eigen space. Now, we
have $dim(D)=2n$ and $dim(H)$ is even, which shows that $dim(H\cap
D)\geq 1$. Let $X\in H\cap D$ is a unit vector, then $\phi X$ is
also unit vector then from equation (3.1) we obtain
$$\cos\theta=\sqrt{-\lambda(X)},$$
which proves the first part. Moreover, for any $X\in H$, formula (3.1) yields $\| \phi X \|= \|X\|$ which imply that
$g(X,\xi)=0 $, hence $H\subset D.~~~~~~~~~~~~~~~~~~~~~~~~~~~~~~~~~~~~~~~~~~~~~~~~~~~~~\square$

\parindent=8mm
We have noted that, invariant and anti-invariant submanifolds are slant submanifolds
with slant angle $\theta =0 $ and $\theta= \pi /2$, respectively. A
slant immersion which is neither invariant nor anti-invariant is
called a proper slant immersion.  In case of invariant submanifold
$T=\phi$ and so
$$T^2=\phi^2 = -I+ \eta\otimes\xi.$$
While in case of anti-invariant submanifold, $T^2=0$. In fact, we have the following general result which characterize slant immersion.\\

\noindent{\bf{Theorem 3.3.}} {\it{Let $M$ be a space like submanifold of a Lorentzian manifold $\bar M$ such that
$\xi\in TM$. Then, $M$ is slant submanifold if and only if there
exist a constant $\lambda \in [0,1]$ such that
$$ T^2=\lambda(-I+\eta\otimes\xi).\eqno(3.3)$$
Furthermore, if $\theta $ is slant angle of $M$, then $\lambda =
\cos^2\theta$.}}\\

\noindent{\it{Proof.}} Necessary condition is
obvious, we have to prove the sufficient condition, suppose that
there exist a constant $\lambda$ such that $T^2= \lambda(-I+\eta\otimes \xi)$, then for any $X\in TM-\langle\xi\rangle$, we have
$$\cos\theta(X)=\frac{g(\phi X, TX)}{\|TX\|\|\phi X\|}~~~~~~~~~~$$
$$~~~~~~~~~~~~~~~~~~~=-\frac {g(X,T^2X)}{{\|TX\|\|\phi X\|}}=\lambda\frac{\|\phi X\|}{\| TX\|}.\eqno(3.4)$$
On the other hand $\cos\theta(X)=\frac{\|TX \|}{\|
\phi X\|}$ and so by using (3.4) we obtain that $\lambda=\cos^2\theta $. Hence, $\theta(X)$ is a constant angle of $M$ i.e, $M$ is slant.$~~~~~~~~~~~~~~~~~~~~~~~~~~~~~~~~~~~~~~~~~~~\square$

\parindent=8mm Now, we have the following corollary, which can be easily
verified.\\

\noindent{\bf{Corollary 3.1.}} {\it{Let $M$ be a space
like slant submanifold of a Lorentzian manifold $\bar M$ with slant
angle $\theta$. Then for any $X,Y \in TM$, we have}}
$$ g(TX,TY) = \cos^2 \theta (g(X,Y) - \eta(X) \eta (Y))\eqno(3.5)$$
$$ g(NX,NY) = \sin^2 \theta (g(X,Y) - \eta(X) \eta (Y)).\eqno(3.6)$$

\noindent{\it{Proof.}} For any $X, Y\in TM$, then by equation (2.3) we have
$$g(X, TY)=-g(TX, Y).$$
Substituting $Y$ by $TY$ in the above equation we get
$$g(TX, TY)=-g(X, T^2Y).$$
Then by virtue of (3.3), we obtain (3.5). The proof of (3.6) follows from (2.2) and (2.9).$~~~~~~~~~~~~~~~~~~~~~~~~~~~~~~~~~~~~~~~~~~~~~~~~~~~~~~~~~~~~~~~~~~~~~~~~~~~~~~~~~~~~~~~~~~~~~~~~~~~~\square$

\section{Slant submanifolds of Lorentzian Sasakian manifolds}

\noindent
In this section we study the slant submanifold of Lorentzian Sasakian manifolds and obtain some interesting results using the equation of curvature tensor.\\

\noindent{\bf{Theorem 4.1.}} {\it{Let $M$ be a slant
 submanifold of a  Lorentzian Sasakian manifold $\bar M$. Then $Q$ is parallel
 if and only if $M$ is anti-invariant.}}\\

\noindent{\it{Proof.}} Let $\theta$ be slant angle of $M$ in $\bar
M$, then for any $X, Y\in TM$ by equation (3.3), we have
$$T^2 Y=QY=\cos^2 \theta(-Y+\eta(Y)\xi)\eqno(4.1)$$
and
$$Q\nabla_XY=\cos^2\theta(-\nabla_X Y+ \eta(\nabla_XY)\xi).\eqno(4.2)$$
Taking the covariant derivative of (4.1) with respect to $X\in TM$, we get
$$\nabla_XQY=\cos^2\theta(-\nabla_XY+ \eta(\nabla_XY)\xi+g(Y, \nabla_X\xi)\xi+\eta(Y)\nabla_X\xi).\eqno(4.3)$$
Then from equations (4.2) and (4.3) we get
$$(\bar\nabla_XQ)Y=\cos^2\theta(g(Y, \nabla_X\xi)\xi+\eta(Y)\nabla_X \xi).\eqno(4.4)$$
Thus on using (2.5), (2.6) and (2.9) the above equation takes the form
$$(\bar\nabla_XQ)Y=\cos^2\theta(g(Y, TX)\xi+\eta(Y)TX).\eqno(4.5)$$
The assertion follows from (4.5).$~~~~~~~~~~~~~~~~~~~~~~~~~~~~~~~~~~~~~~~~~~~~~~~~~~~~~~~~~~~~~~~~~~~~~~~~~~~~~\square$

\parindent=8mm Now, we shall investigate the existence of a slant submanifold using curvature tensor.\\

\noindent{\bf{Lemma 4.1.}} {\it{Let $M$ be a submanifold of Lorentzian Sasakian manifold $\bar M$ such that $\xi$ is tangent to $M$. Then for any $X, Y\in TM$, we have
$$R(X, Y)\xi=(\nabla_Y T)X-(\nabla_XT)Y\eqno(4.6)$$
where $R$ is the curvature tensor field associated to the metric induced by $\bar M$
on $M$. Moreover,}}
$$R(\xi, X)\xi=QX-(\nabla_\xi T)X\eqno(4.7)$$
$$R(X,\xi, X, \xi)=g(QX, X).\eqno(4.8)$$

\noindent{\it{Proof.}} For any $X\in TM$ then from (2.5) and (2.9) we have
$$TX=\nabla_X \xi.$$
Using this fact in the formula $(\nabla_XT)Y=\nabla_XTY-T\nabla_XY$, we obtain
$$(\nabla_XT)Y=\nabla_X\nabla_Y\xi-\nabla_{\nabla_XY}\xi.$$
Similarly,
$$(\nabla_YT)X=\nabla_YTX-T\nabla_YX=\nabla_Y\nabla_X\xi-\nabla_{\nabla_Y X}\xi.$$
Substituting these equations in the definition of $R(X,Y)\xi$ it is easy to get (4.6). Rewriting (4.6) for $X=\xi$ and $Y=X$, we obtain
$$R(\xi, X)\xi=(\nabla_X T)\xi-(\nabla_\xi T)X= QX-(\nabla_\xi T)X.$$
Which proves (4.7). Now taking the product with $X$ in (4.7), we get
$$R(\xi, X, \xi, X)=g(QX, X)-g((\nabla_\xi T)X, X).\eqno(4.9)$$
The second term of (4.9) will be identically zero as follows
$$g((\nabla_\xi T)X, X)=g(\nabla_\xi TX, X)-g(T\nabla_\xi X, X)~~~~~~~~~~~~~~~~$$
$$~~~~~~~~~~~~=-g(TX,\nabla_\xi X)+g(\nabla_\xi X, TX)=0.$$
Then (4.8) follows from (4.9) using the above fact.$~~~~~~~~~~~~~~~~~~~~~~~~~~~~~~~~~~~~~~~~~~~~~~~~~~~~~~\square$\\

\noindent{\bf{Theorem 4.2.}} {\it{Let $M$ be a submanifold of a Lorentzian Sasakian manifold $\bar M$ such that the
characteristic vector field $\xi $ is tangent to $M$. If $\theta \in(0, \pi/2)$ then the following statements are equivalent}}
\begin{enumerate}
\item [{(i)}] {\it{$M$ is slant with slant angle $\theta$.}}
\item [{(ii)}] {\it{For any $x\in M$ the sectional curvature of any $2$-plane of
$T_x M$ containing $\xi_x$ equals $\cos^2\theta$.}}
\end{enumerate}

\noindent{\it{Proof.}} Assume that the statement $(i)$ is true, then for any $X
\perp \xi$ by Theorem 3.3, we have
$$QX=\cos^2\theta X$$
which by virtue of (4.8) yields
$$R(X,\xi, X, \xi)=\cos^2\theta.\eqno(4.10)$$
Thus $(ii)$ is proved.

\par Conversely, suppose that $(ii)$ hold then for any $X\in TM$, we may write
$$X=X_\xi+{X_\xi}^\perp \eqno (4.11) $$
where $X_\xi=\eta(X)\xi$ and ${X_\xi}^\perp$ is the
component of $X$ perpendicular to the $\xi$, using (4.10) and (4.11)
$$\frac{R({X_\xi}^\perp, \xi, {X_\xi}^\perp, \xi)}{\vert{{X_\xi}\vert}^2}=\cos^2 \theta,$$
$${\mbox{or}},~~~~~~~~~~~~~~~~~~~~~~~~~~~~~~~~~~R({X_\xi}^\perp,\xi,{X_\xi}^\perp,\xi)=\cos^2 \theta {\vert X_\xi \vert}^2.~~~~~~~~~~~~~~~~~~~~~~~~~~~~~~~~~~~~~~~~~~~~~\eqno(4.12)$$
Let $X$ be a unit vector filed such that $QX=0$. Then from (4.8) and (4.12)
$$\cos^2\theta\vert {X_\xi}^\perp\vert^2=0.\eqno(4.13)$$
If $\cos\theta\neq0$, then from the above equation $X=X_\xi.$ This proves that at each point $x\in M$,
$$Ker(Q)=\langle \xi_x \rangle .\eqno(4.14)$$
Moreover, Let $A$ be the matrix of the endomorphism
$Q$ at $x\in M$, then for a unit vector field $X$ on $M$, $QX=AX,$
and as $Q(X_\xi)=0, X=X_\xi.$ Then by (4.8) and (4.12)
$$ A=\cos^2\theta I.\eqno(4.15)$$
Choosing $\lambda=\cos^2\theta$, we conclude that for any
$x\in M$, This fact together with (4.14) and Theorem 3.3, verifies
that $M$ is slant in $\bar M$ with slant angle $\theta$. Finally, suppose $\cos\theta=0$ and $X$ is an arbitrary unit vector field such that $QX=\lambda X$ where $\lambda\in C^\infty(M)$. Then, from (4.8) and (4.12) $g(QX, X)=0$ that is $\lambda=0$ and therefore $Q=0$ which means that $M$ is anti-invariant.$~~~~~~~~~~~~~~~~~~~~~~~~~~~~~~~~~~~~~~~~~~~~~~~~~\square$

\bigskip
Khushwant Singh, Satvinder Singh Bhatia:

\noindent School of Mathematics and Computer Applications\\
Thapar University, Patiala-147 004, INDIA\\
\noindent {\it E-mail}: {\tt khushwantchahil@gmail.com,
ssbhatia@thapar.edu}

\bigskip

Siraj Uddin:

\noindent Institute of Mathematical Sciences,\\
Faculty of Science,University of Malaya,\\
50603 Kuala Lumpur, MALAYSIA

\noindent {\it E-mail}: {\tt siraj.ch@gmail.com}

\bigskip

Meraj Ali Khan:

\noindent Department of Mathematics,\\
University of Tabuk, Tabuk,\\
Kingdom of Saudi Arabia.

\noindent {\it E-mail}: {\tt meraj79@gmail.com}

\smallskip
\end{document}